\documentclass[10pt]{amsart}

\usepackage{amsmath, amscd, amssymb}
\numberwithin{equation}{section}
\newcommand{\sM}{{\mathcal M}}

\newtheorem{theorem}{Theorem}[section]

\newtheorem{proposition}[theorem]{Proposition}
\newtheorem{lemma}[theorem]{Lemma}

\newtheorem{corollary}[theorem]{Corollary}

\newtheorem{conjecture}{Conjecture}

\theoremstyle{remark}

\theoremstyle{definition}

\begin{document}

\title[Bouchard-Klemm-Marino-Pasquetti Conjecture for $\mathbb{C}^3$]
{Bouchard-Klemm-Marino-Pasquetti Conjecture for $\mathbb{C}^3$}

\author{Lin Chen}
\address{Simons Center for Geometry and Physics\\
State University of New York, Stony Brook}
\email{chenlin@math.sunysb.edu}

\begin{abstract}
In this paper, we give a proof of the
Bouchard-Klemm-Marino-Pasquetti conjecture for a framed vertex, by
using the symmetrized Cut-Join Equation developed in a previous
paper.

\end{abstract}

\maketitle

\setcounter{tocdepth}{5} \setcounter{page}{1}

\section{Introduction}

In their seminal paper \cite{BKMP}, Bouchard, Klemm, Marino and
Pasquetti propose a new approach to compute both the open and closed
Gromov-Witten invariants of local Calabi-Yau manifolds, including
the mirrors of toric varieties. The approach is based on the matrix
models theory of Eynard and Orantin. To each toric Calabi-Yau
three-fold, there is an algebraic curve $\Sigma$, living inside
$\mathbb{C}^{*}\times\mathbb{C}^{*}$, associated to the toric
diagram. The curve $\Sigma$ is called the framed mirror curve, with
genus equal to the genus of the toric diagram. Bouchard, Klemm,
Marino and Pasqquetti conjecture that the Gromov-Witten invariants
of the toric Calabi-Yau three folds can be computed by applying the
recursion of Eynard and Oranin to the framed mirror curve.

The first example of Bouchard-Klemm-Marino-Pasquetti theory is the
non-compact toric three fold $\mathbb{C}^3$, the socalled framed
vertex, which is the building block of three dimensional toric
varieties. In this case, the framed mirror curve is the algebraic
curve $x=y^f(1-y)$. The theory of Eynard and Orantin produce a
topological recursion relations, when plug in the information of the
framed mirror curve. On the other hand, there is an existing
topological recursion for $\mathbb{C}^3$, i.e, the cut-join equation
proved in \cite{LLZ1} and its symmetrized version \cite{C}. In this
paper, we will prove that the Symmetrized Cut-Join Equation obtained
in \cite{C} implies the topological recursion of Eynard and Orantin,
thus prove the Bouchard-Klemm-Marino-Pasquetti conjecture in the
$\mathbb{C}^3$ case.

Our strategy is strongly motivated by the recent work of
Eynard-Mulase-Safnuk on the Bouchard-Marino conjecure. The
Bouchard-Marino conjecture on Hurwitz numbers is a topological
recursion of Eynard-Orantin type. This conjecture was first proved
by Borot, Eynard, Mulase and Safnuk, using very deep results in
matrix model theory. Later, Eynar, Mulase and Safnuk give another
proof by comparing the conjectural recursion with the symmetrized
cut-join equation of Hurwitz numbers, which was discovered by
Goulden-Jackson-Vainshtein \cite{GJV}. We now describe their second
approach. First, writing the residues in the conjectural recursion
as a contour integral, and a residue theorem calculus switch the
calculation to the two nearby simple poles. This computation provide
an equivalent form of the conjectural formula, which is an identity
of polynomials. Then, by pushing forward the Symmetrized Cut-Join
Equation of Goulden-Jackson-Vainstein via the projection
$\pi:\Sigma'\rightarrow\mathbb{C}$ from the mirror curve $\Sigma'$
associated to the Hurwitz numbers, they are able to show the
resulting equation, modulo the principle (singular) part, is
precisely the equation obtained in the first step.

In this paper, we will give a proof of the BKMP conjecture for the
framed vertex $\mathbb{C}^3$. In fact, the Bouchard-Marino
conjecture is a specialization of the BKMP conjecture for
$\mathbb{C}^3$, by letting the framing $f\rightarrow \infty$. The
Hurwitz numbers are then replaced by the Gromov-Witten invariants,
which can be written as Hodge integrals involves three $\lambda$
classes. The generating series of such Hodge integrals satisfies a
similar cut-join equation. The corresponding Symmetrized Cut-Join
Equation is obtained by the author in a previous paper \cite{C}.
Following the line of Eynard-Mulase-Safnuk, we prove the BKMP
conjecture by switching the residues calculation to the two nearby
simple poles in the question, and pushing forward the Symmetrized
Cut-Join Equation of \cite{C} via the projection from the framed
mirror curve.

We will describe some known results in Section \ref{Sec2}. The BKMP
conjecture for $\mathbb{C}^3$ is stated in Section \ref{Sec3}. In
Section \ref{Sec4}, we compute the residues appeared in the
conjectural recursion. Replacing the complex analysis of the
functions $\eta_n(v)$, invent by Eynard-Mulase-Safnuk, by formal
power series argument in Section \ref{Sec5}, we simplify the proof
the technical Lemma \ref{Lemma3}, avoid the convergence and analytic
continuation difficulty in \cite{EMS}. Finally, we compute the push
forward of the Symmetrized Cut-Join Equation in the last two
sections, and establish the BKMP conjecture.

I would like to thank Professor Kefeng Liu and Bailin Song for
helpful discussion. I also want to thank Professor Mulase for
explaining their work to me. This work is supported by Simons Center
for Geometry and Physics, Stony Brook University.

\vskip 1cm

\section{Marino-Vafa formula and Symmetrized Cut-Join Equation}
\label{Sec2}

The celebrated topological vertex theory developed by
Li-Liu-Liu-Zhou \cite{LLLZ} establishes a correspondence between two
different types of physics theory, topological string theory and the
Chern-Simons theory. More precisely, it gives a closed formula of
Gromov-Witten invariants of certain toric Calabi-Yau three folds in
terms of link invariants (see \cite{LLZ1} \cite{LLZ2} \cite{LLLZ}).

The simplest example is the so-called Marino-Vafa formula proved by
Liu-Liu-Zhou \cite{LLZ1}. Via virtual localization, one side of the
Marino-Vafa formula can be written as a generating series of Hodge
integrals involves three lambda classes
$\mathcal{C}=\sum_{g\geqslant 0,n\geqslant
1}\mathcal{C}_n^g\lambda^{2g-2+n}$, where
\begin{align*}
\mathcal{C}_n^g=&\sum_{d\geqslant 1}\sum_{\mu\vdash
d,l(\mu)=n}-\frac{\sqrt{-1}^{d+n}}{|\mathrm{Aut}
\mu|}(f(1+f))^{n-1}\prod_{i=1}^n\frac{\prod_{a=1}^{\mu_i-1}
(\mu_if+a)}{(\mu_i-1)!}\int_{\overline{\sM}_{g,n}}\frac{\Gamma_g(f)}
{\prod_{i=1}^n(1-\mu_i\psi_i)}\cdot \mathbf{p}_\mu\\
=&-\frac{\sqrt{-1}^{n}(f(1+f))^{n-1}}{n!}\sum_{\mu_1,\mu_2\cdots\mu_n\geqslant
1}^{+\infty}\sqrt{-1}^{|\mu|}\prod_{i=1}^n\frac{\prod_{a=1}^{\mu_i-1}(\mu_if+a)}{(\mu_i-1)!}\\
&\cdot\sum_{b_1,\cdots,b_n}\int_{\overline{\sM}_{g,n}}
\Gamma_g(f)\prod_{i=1}^n\psi_i^{b_i}\prod_{i=1}^n\mu_i^{b_i}\cdot \mathbf{p}_\mu\\
=&-\frac{(f(1+f))^{n-1}}{n!}\sum_{b_1,\cdots,b_n}
<\tau_{b_1}\cdots\tau_{b_n}\Gamma_g(f)>\prod_{i=1}^n\varphi_{b_i}(\overrightarrow{\mathbf{p}}).
\end{align*}
The formula proved in \cite{LLZ1} relates the generating series
$\mathcal{C}$ to a truncated version of the framing dependent
Chern-Simons partition function (generating series of colored HOMFLY
polynomials) of unknot. In the above formula, we use the notation
$$
\varphi_i(\overrightarrow{\mathbf{p}})=\sum_{m\geqslant
1}\sqrt{-1}^{m+1}
\mathbf{p}_m\frac{\prod_{a=1}^{m-1}(mf+a)}{(m-1)!}m^i=
\frac{1}{f}\sum_{m\geqslant 1}\sqrt{-1}^{m+1}\mathbf{p}_m
\frac{\prod_{a=0}^{m-1}(mf+a)}{m!}m^i
$$
for a formal sum involve infinitely many formal variables
$\overrightarrow{\mathbf{p}}=\{\mathbf{p}_1,\mathbf{p}_2,\cdots\}$,
and the class
$$
\Gamma_g(f)=\Lambda_g^\vee(1)\Lambda_g^\vee(f)\Lambda_g^\vee(-f-1)
$$
for
$$
\Lambda_g^\vee(u)=u^g-\lambda_1 u^{g-1}+\cdots+(-1)^g\lambda_g.
$$
The relation between these Hodge integrals and the Gromov-Witten
invariants of Calabi-Yau three-folds can be found in \cite{KLiu}.

The generating function $\mathcal{C}$ satisfy a Cut-Join Equation:
\begin{equation} \label{E:CJ1}
\frac{\partial\mathcal{C}}{\partial
f}=\frac{\sqrt{-1}\lambda}{2}\sum_{i,j\geqslant
1}(ij\mathbf{p}_{i+j}\frac{\partial^2\mathcal{C}}
{\partial\mathbf{p_i}\partial\mathbf{p_j}}+ij\mathbf{p}_{i+j}
\frac{\partial\mathcal{C}}{\partial\mathbf{p}_i}
\frac{\partial\mathcal{C}}{\partial\mathbf{p}_j}+(i+j)
\mathbf{p}_i\mathbf{p}_j\frac{\partial\mathcal{C}}{\partial\mathbf{p}_{i+j}})
\end{equation}

In \cite{C}, the author introduced symmetrization operators
$$
\Xi_n:\mathbb{C}[\mathbf{p}_1,\mathbf{p}_2,\cdots]\rightarrow\mathbb{C}[x_1,\cdots,x_n]
$$
with values in the subring of $\mathbb{C}[x_1,\cdots,x_n]$ consists
of symmetric polynomials, by letting
$$
\Xi_n \mathbf{p}_\alpha=(\sqrt{-1})^{-n-|\alpha|}\sum_{\sigma\in
S_n}x_{\sigma(1)}^{\alpha_1}\cdots x_{\sigma(n)}^{\alpha_n}
$$
for $n\geqslant 1$ if $l(\alpha)=n$ with
$\alpha=(\alpha_1,\cdots,\alpha_n)$, and 0 otherwise.

\vskip 0.5cm

After a transcendental change of variables
\begin{align*}
&y(x)^f(1-y(x))=x\\
&t=\frac{1}{(1+f)y-f},
\end{align*}
or explicitly
\begin{align*}
t&=1+(\frac{1+f}{f})((1+f)x\frac{d}{dx}(1-y)-(1-y))\\
&=1+(\frac{1+f}{f})(\frac{f(1-y)}{(1+f)y-f})\\
&=1+(\frac{1+f}{f})\sum_{n=1}^{+\infty}\frac{\prod_{a=0}^{n-1}(nf+a)}{n!}x^n,
\end{align*}
the generating series become
$$
H_{g,n}= -{(f(1+f))^{n-1}}\sum_{b_1,\cdots,b_n}
<\tau_{b_1}\cdots\tau_{b_n}\Gamma_g(f)>_g\prod_{i=1}^n\phi_{b_i}(t_i).
$$
In the above expression,
$$
\phi_b(t)=\left(\frac{t(t-1)(ft+1)}{f+1}\frac{d}{dt}\right)^b\left(\frac{t-1}{f+1}\right)
$$
for $b\geq 0$ is a polynomial in $t$ of degree $2b+1$, with
coefficients in $\mathbb{Q}(f)$. For convenience, we denote
$$
\phi_{-1}(t)=-\log(1+\frac{1}{ft})
$$
for later use. It is easy to check that
$$
\phi_0(t)=\frac{t-1}{f+1}=\left(\frac{t(t-1)(ft+1)}{f+1}\frac{d}{dt}\right)\phi_{-1}(t)
$$
compatible with our notations.

\vskip 0.5cm

Apply the symmetrization operator $\Xi_n$ and change the variables
from $x_i$ into $t_i$, Equation \ref{E:CJ1} transformed into the
following Symmetrized Cut-Join Equation proved in \cite{C}
\begin{align*}
(\frac{\partial}{\partial
f}&+\sum_{l=1}^n\frac{t_l(t_l-1)}{f+1}\cdot \frac{\partial}{\partial
t_l})H_n^g(t_1,\cdots,t_n,f)=T_1+T_2+T_3+T_4,
\end{align*}
where
\begin{align*}
T_1=&-\frac{1}{2}\sum_{l=1}^nt_l(t_l-1)(\frac{t_lf+1}{f+1})\frac{\partial}{\partial
t_l}\cdot
t_{n+1}(t_{n+1}-1)(\frac{t_{n+1}f+1}{f+1})\frac{\partial}{\partial
t_{n+1}}H_{n+1}^{g-1}|_{t_{n+1}=t_l}\\
T_2=&-\frac{1}{2}\sum_{1\leqslant a\leqslant g-1}\sum_{1\leqslant
k\leqslant
n}\Theta_{k-1}(t_1(t_1-1)(\frac{t_1f+1}{f+1})\frac{\partial}{\partial
t_1}H_k^a(t_1,\cdots,t_k,f))\\
&\cdot(t_1(t_1-1)(\frac{t_1f+1}{f+1})\frac{\partial}{\partial
t_1}H_{n-k+1}^{g-a}(t_1,t_{k+1},\cdots,t_n,f))\\
T_3=&-\sum_{k=3}^n\Theta_{k-1}(t_1(t_1-1)(\frac{t_1f+1}{f+1})\frac{\partial}{\partial
t_1}H_k^0(t_1,\cdots,t_k,f))(t_1(t_1-1)
(\frac{t_1f+1}{f+1})\frac{\partial}{\partial
t_1}H_{n-k+1}^{g}(t_1,t_{k+1},\cdots,t_n,f))\\
T_4=&\Theta_1\frac{t_1^2(t_1-1)(t_2-1)}{t_1-t_2}(\frac{t_1f+1}{f+1})^2\frac{\partial}{\partial
t_1}H_{n-1}^{g}(t_1,t_3,\cdots,t_n,f).
\end{align*}
We remark that this is an equality in the polynomial ring
$\mathbb{C}(f)[t_1,\cdots,t_n]$. The symbol $\Theta_k$ means take
the symmetric sum. For more detail about this Symmetrized Cut-Join
Equation, we refer to \cite{C}.

\vskip 1cm

\section{The BKMP Conjecture} \label{Sec3}

The Symmetrized Cut-Join Equation is a set of topological recursion
of Hodge integrals, i.e, it computes genus $g$ Hodge integrals in
terms of integrals with less genus or less marked points.

The BKMP conjecture is another set of topological recursion of Hodge
integrals, coming from the conjectural duality between topological
string theory and matrix model theory. We now describe the
recursions in the BKMP conjecture for $\mathbb{C}^3$. In this case,
the Framed Mirror Curve $\mathrm{C}$ is given by
$$
x=y^f(1-y),
$$
with the formal power series inverse function
$$
y(x)=1-\sum_{n=1}^{\infty}\frac{\prod_{j=0}^{n-2}(nf+j)}{n!}x^n.
$$
Denote $\zeta_n$ (depends on $f$) the differential one form
\begin{align*}
\zeta_n&=d(\phi_n)=d\left((x\frac{d}{dx})(\frac{t-1}{f+1})\right)\\
&=d\left((\frac{t(t-1)(ft+1)}{f+1})^n(\frac{t-1}{f+1})\right)\\
&=d(\frac{1}{\tau}\sum_{m=1}^{\infty}\frac{\prod_{a=0}^{m-1}(mf+a)}{m!}m^nx^m).
\end{align*}

Let $W_g(x_1,\cdots,x_n)$ be the differential $n$-form defined by
the formula
\begin{align*}
W_g(x_1,\cdots,x_n)=(-1)^{g+n}(f(f+1))^{n-1}
\sum_{b_1,\cdots,b_n}<\tau_{b_1}\cdots\tau_{b_n}\Gamma_g(f)>
\prod_{i=1}^n\zeta_{b_i}(y_i(x_i),f).
\end{align*}
For example
\begin{align*}
&W_0(x)=\log y(x)\cdot \frac{dx}{x}\\
&W_0(x_1,x_2)=B(y_1,y_2)-\frac{dx_1dx_2}{(x_1-x_2)^2}
=\frac{dy_1dy_2}{(y_1-y_2)^2}-\frac{dx_1dx_2}{(x_1-x_2)^2}.\\
&W_0(x_1,x_2,x_3)=-\frac{f^2}{(f+1)}dt_1dt_2dt_3
\end{align*}
If we denote $d_i=\frac{\partial}{\partial
x_i}dx_i=\frac{\partial}{\partial t_i}dt_i$ the differential with
respect to the $i$-th variable, then the differential $n$-form
$$
W_g(x_1,\cdots,x_n)=(-1)^{g+n-1}d_1\cdots d_n H_g^n
$$
becomes an element in the canonical module
$\mathbb{C}(f)[t_1,\cdots,t_n]dt_1\cdots dt_n$, if switch to the
$t_i$ variables.

\vskip 0.5cm

The Framed Mirror Curve $\mathrm{C}$ has a critical point
$(x,y)=(\frac{f^f}{(f+1)^{f+1}},\frac{f}{f+1})$. Near the critical
point, the morphism
$$
\mathrm{C}\rightarrow\mathbb{C}
$$
sending $(x,y)$ to $x$ is (locally) a branching double cover. Let
$q$ and $\overline{q}$ be two points on the Framed Mirror Curve
close to the critical point such that $x(q)=x(\overline{q})$. Define
one forms
\begin{align*}
&\omega(q)=\left(\log y(q)-\log
y(\overline{q})\right)\frac{dx(q)}{x(q)}\\
&dE(q,\overline{q},y_2)=\frac{dy_2}{2}
\left(\frac{1}{y_2-y(q)}-\frac{1}{y_2-y(\overline{q})}\right).
\end{align*}
The kernel function is defined to be the formal quotient
$$
K(q,\overline{q},y_2)=\frac{dE(q,\overline{q},y_2)}{\omega(q)}.
$$

The BKMP conjectural recursion reads:
\begin{conjecture}[BKMP] \label{Conj1}
The differential forms $W_g(y,y_1,\cdots,y_n)$ are completely
determined by the following topological recursion relation
\begin{align*}
W_g(y,y_1,\cdots,y_n)=&\mathrm{Res}_{y(q)=\frac{f}{f+1}}\frac{dE(q,\overline{q},y)}{w(q)}
[W_{g-1}(y(q),y(\overline{q}),y_1,\cdots,y_n)\\
&+\sum_{g_1+g_2=g}\sum_{I\coprod
J=H}^{\mathrm{stable}}\sum_{J\subset
H}W_{g_1}(y(q),y_I)W_{g_2}(y(\overline{q}),y_{J})\\
&+\sum_{i=1}^n (W_g(y(q),y_{H\setminus\{i\}})\otimes
B(y(\overline{q}),y_i)+B(y(q),y_i)\otimes
W_g(y(\overline{q}),y_{H\setminus\{i\}}))],
\end{align*}
together with the initial conditions
\begin{align*}
&W_0(x_1,x_2,x_3)=-\frac{f^2}{(f+1)}dt_1dt_2dt_3\\
&W_1(x_1)=\frac{1}{24}((1+f+f^2)\zeta_0-f(1+f)\zeta_1).
\end{align*}
\end{conjecture}

In the rest of this paper, we will prove that the Symmetrized
Cut-Join Equation implies the above BKMP conjecture.

\vskip 1cm

\section{Residue Calculus} \label{Sec4}

The RHS of the BKMP conjecture involves the following two types of
residues.

Type I: Comes from $W_{g-1}$ and the stable sums.
$$
R_{a,b}(y)=\mathrm{Res}_{y(q)=\frac{f}{f+1}}
\frac{dE(q,\overline{q},y)}{\omega(q)}\zeta_a(y(q))\zeta_b(y(\overline{q})).
$$

Type II: Comes from the unstable contribution of integral over the
point $\overline{M}_{0,2}$.
$$
R_{a}(y,y_i)=\mathrm{Res}_{y(q)=\frac{f}{f+1}}\frac{dE(q,\overline{q},y)}{\omega(q)}
\left(\zeta_a(y(q))B(y(\overline{q}),y_i)+\zeta_a(y(\overline{q}))B(y(q),y_i)\right).
$$

\vskip 0.5cm

Before going into the computations, we first explain the meaning of
the above two residues. Take a small open neighborhood $U$ of
$\frac{f^f}{(f+1)^{f+1}}$ in $\mathbb{C}$, such that the projection
$$
\pi:\mathrm{C}\rightarrow\mathbb{C}
$$
is two to one on the open set
$\pi^{-1}(U\setminus\{\frac{f^f}{(f+1)^{f+1}}\})$. Denote by
$s:q\mapsto \overline{q}$ the holomorphic inversion, which is well
defined on the open neighborhood $\pi^{-1}U$ of the critical point.

In the formula of Type I case,
\begin{align*}
&\zeta_a(y(q))=\phi_a'(y(q))dx(q)\\
&\zeta_b(y(\overline{q}))=\phi_b'(y(\overline{q}))s'(x(q))dx(q)
\end{align*}
where $\phi_n'=\frac{d}{dx}\phi_n$ is the derivative of $\phi_n$.
The denominator $w(q)$ of the formal quotient
$K(q,\overline{q},y_2)$ canceled one of the $dx(q)$, thus the
product $K(q,\overline{q},y)\zeta_a(y(q))\zeta_b(y(\overline{q}))$
is a two form. After taking the residue, $R_{a,b}(y)$ became a
meromorphic one form. Similar description holds for the two form
$R_a(y,y_i)$.

\vskip 0.8cm

\begin{lemma} \label{Lemma1}
Change into the variable $t$, i.e, let $R_{a,b}(y)=P_{a,b}(t)dt$,
then $P_{a,b}(t)\in\mathbb{C}(f)[t]$ is a polynomial of $t$ with
coefficients in the field $\mathbb{C}(f)$.
\end{lemma}


\begin{proof}
Let $y=\frac{f+\frac{1}{t}}{f+1}$ and let $z=\frac{1}{t}$. The
critical point of the Framed Mirror Curve is then the point $0$ in
the $z$ coordinate. Formally, the kernel function
$K(q,\overline{q},y)$ is the product of
$\frac{(f+1)dy}{2}\otimes\frac{x(q)}{dx(q)}$ and the function
\begin{align*}
\frac{(\frac{1}{z-z(q)}-\frac{1}{z-z(\overline{q})})}
{\log(1+\frac{z(q)}{f})-\log(1+\frac{z(q)}{f})}=\frac{1}{(z-z(q))(z-z(\overline{q}))}
\cdot
[\frac{z(q)-z(\overline{q})}{\log(1+\frac{z(q)}{f})-\log(1+\frac{z(q)}{f})}]
\end{align*}
holomorphic around $z(q)=0$.

We abuse the notation
$$
\zeta_a(y(q))=\zeta_a(t(q))=\zeta_a(\frac{1}{z(q)}),
$$
which is in the module of differentials $\mathbb{C}(f)[t(q)]dt(q)$.

Recall that
$$
\frac{x(q)}{dx(q)}=\frac{t(q)(t(q)-1)(ft(q)+1)}{(f+1)dt(q)}=-\frac{(1-z)(f+z)}{z(f+1)dz}.
$$
By linearity, we only need to consider the residue
$$
\mathrm{Res}_{z(q)=0}\frac{F(z(q))}{(z-z(q))(z-z(\overline{q}))}\cdot
\frac{1}{z(q)^a}\frac{dz(q)}{z(q)^5}\otimes dy
$$
for integer $a\geq 0$ and functions $F(z)$ holomorphic around $z=0$.
Since $dy=-\frac{dt}{t^2(f+1)}$, we have
\begin{align*}
-\frac{dt}{f+1}\cdot\mathrm{Res}_{z(q)=0}\left(F(z(q))
\cdot\sum_{k=0}^{+\infty}\frac{1}{z^{k}}(\sum_{i=0}^kz(q)^iz(\overline{q})^{k-i})
\frac{dz(q)}{z(q)^{5+a}}\right)
\end{align*}
by expanding the denominator. This expression only contains
non-negative power of $\frac{1}{z}=t$. Moreover, if $k\geq 5+a$,
then the summand is then holomorphic at $z(q)=0$, and has no
residue. This shows that the above residue is a polynomial times
$dt$. A more careful analysis of the parameter $f$ in the above
process shows that the coefficients are in fact rational functions
in the variable $f$.
\end{proof}

\vskip 0.4cm






Under the $z$ coordinate, the projection
$\pi:\mathrm{C}\rightarrow\mathbb{C}$ maps the critical point to
$z=0$. By shrinking the open neighborhood $U$ suitably, we can find
a simple closed curve around $z=0$. The holomorphic inversion $s$ is
thus defined in the interior of the area bounded by $\gamma$. Let
the variable $y$ sufficiently close to the critical point such that
it is inside the curve $\gamma$.

Let $\zeta_a(y)=F_a(t)dt$ for some polynomial
$F_a(t)\in\mathbb{C}(f)[t]$. By compactness, the functions
$F_a(y(q))$, $F_b(y(\overline{q}))$,
$\frac{1}{\log(y(q))-\log(y(\overline{q}))}$ and
$t(q)(t(q)-1)(ft(q)+1)$ are all bounded along the circle $\gamma$.
We thus have the following estimate of the contour integral:
\begin{align*}
&\frac{1}{2\pi
i}\int_{\gamma}\frac{dE(q,\overline{q},y)}{\omega(q)}\zeta_a(y(q))\zeta_b(y(\overline{q}))\\
&=\frac{1}{2\pi
i}\int_{\gamma}\frac{\frac{1}{2}(\frac{1}{y-y(q)}-\frac{1}{y-y(\overline{q})})}{\log
y(q)-\log y(\overline{q})}\cdot (\frac{x(q)}{dx(q)}\otimes dy)\cdot
F_a(t(q))dt(q) \cdot F_b(t(\overline{q}))dt(\overline{q})\\
&=\frac{dy}{2\pi
i}\int_{\gamma}\frac{\frac{1}{2}(\frac{1}{y-y(q)}-\frac{1}{y-y(\overline{q})})}{\log
y(q)-\log y(\overline{q})}\cdot\frac{t(q)(t(q)-1)(ft(q)+1)}{(f+1)}
\cdot F_a(t(q))\cdot F_b(t(\overline{q}))s'(t(q))dt(q)\\
& \leq M\cdot\frac{dt}{|t^2|}\thicksim\frac{dt}{|t^2|}.
\end{align*}

By hypothesis, we choose $z$ small enough (sufficiently close to the
critical point), such that it is in the interior bounded by
$\gamma$. Then we choose another simple closed curve $\gamma_\infty$
around $z=0$ (i.e, $t=\infty$), such that both $z$ and $s(z)$ are
outside $\gamma_\infty$.

The above estimate together with Lemma \ref{Lemma1} imply that the
above integral is the principle part of the following integral
\begin{align*}
A:=&\frac{1}{2\pi
i}\int_{\gamma\cup-\gamma_{\infty}}\frac{dE(q,\overline{q},y)}{\omega(q)}
\zeta_a(y(q))\zeta_b(y(\overline{q})).\\
\end{align*}
There are two poles $z(q)=z$ and $z(q)=s(z)$ between the two circle
$\gamma$ and $\gamma_{\infty}$. The integral $A$ thus can be
computed by residue theorem:
\begin{align*}
A&=\mathrm{Res}_{z(q)=z,s(z)}
\frac{\frac{f+1}{2}(\frac{1}{z-z(q)}-\frac{1}{z-z(\overline{q})})}
{\log(1+\frac{z(q)}{f})-\log(1+\frac{z(\overline{q})}{f})}
dy\zeta_a(\frac{1}{z(q)})\zeta_b(\frac{1}{z(\overline{q})})\cdot\frac{x(q)}{dx(q)}\\
&=\mathrm{Res}_{z(q)=z,s(z)}
\frac{\frac{z^2dt}{2}(\frac{1}{z-z(q)}-\frac{1}{z-z(\overline{q})})}
{\log(1+\frac{z(q)}{f})-\log(1+\frac{z(\overline{q})}{f})}\cdot
F_a(\frac{1}{z(q)})F_b(\frac{1}{z(\overline{q})})s'(\frac{1}{z(q)})
\frac{dz(q)}{z(q)^2}
\cdot \frac{(1-z(q))(z(q)+f)}{(f+1)z(q)^3}\\
&=\frac{-\frac{1}{2}dt}{\log(1+\frac{1}{ft})-\log(1+\frac{1}{fs(t)})}
F_a(t)F_b(s(t))s'(t)\cdot\frac{t(t-1)(ft+1)}{(f+1)}\\
&+\frac{-\frac{1}{2}dt}{\log(1+\frac{1}{ft})-\log(1+\frac{1}{fs(t)})}
F_a(s(t))F_b(t)\cdot
\frac{s(t)(s(t)-1)(fs(t)+1)}{(f+1)}\\
&=\frac{-\frac{1}{2}}{\log(1+\frac{1}{ft})-\log(1+\frac{1}{fs(t)})}
\left(\phi_{a+1}(t)d\phi_b(s(t))+\phi_{a+1}(s(t))d\phi_{b}(t)\right)\\
&=\frac{-\frac{1}{2}}{\log(1+\frac{1}{ft})-\log(1+\frac{1}{fs(t)})}
\left(\phi_{a+1}(t)\phi_{b+1}(s(t))+\phi_{a+1}(s(t))\phi_{b+1}(t)
\right)\cdot \frac{(f+1)dt}{t(t-1)(ft+1)}.
\end{align*}

\vskip 0.7cm

Let $f(t)=\sum_{n\geq-N}^{\infty}a_nt^n\in\mathbb{C}(f)((t))$ be a
Laurent series. Denote by $f(t)_+=\sum_{n\geq 0}^{\infty}a_nt^n$ the
regular part of $f(t)$. The above formula of $A$ and the estimate of
the contour integral along $\gamma$ give the following proposition.
\begin{proposition} \label{Prop1}
$$
P_{a,b}(t)=\left(\frac{\frac{1}{2}}{\log(1+\frac{1}{ft})-\log(1+\frac{1}{fs(t)})}
\left(\phi_{a+1}(t)\phi_{b+1}(s(t))+\phi_{a+1}(s(t))\phi_{b+1}(t)
\right)\cdot \frac{(f+1)}{t(t-1)(ft+1)}\right)_+
$$
\end{proposition}

\vskip 1cm

Next we consider the Type II residues. A similar argument as in the
Type I case leads to the following lemma, the proof of which we left
for the reader.

\begin{lemma} \label{Lemma2}
Change to the $t$ and $t_i$ variable, $R_a(y,y_i)=P_a(t,t_i)dtdt_i$,
then $P_a(t,t_i)\in\mathbb{C}(f)[t,t_i]$ is a polynomial of two
variables $t$ and $t_i$ with coefficients in the field
$\mathbb{C}(f)$.
\end{lemma}

\vskip 0.5cm

Let the variable $t_i$ take values in the area such that its inverse
$z_i$ and $s(z_i)$ are outside the area bounded by $\gamma$. Let us
now consider the following integral
$$
B:=\frac{1}{2\pi i}\int_{\gamma\bigsqcup-\gamma_{\infty}}
\frac{dE(q,\overline{q},y)}{w(q)}(\zeta_z(y(q))B(y(\overline{q}),y_i)
+\zeta_a(y(\overline{q}))B(y(q),y_i)).
$$

We expect the integral along $\gamma$ is the principle part, while
the $\gamma_{\infty}$ part gives the residue we compute.

By our hypothesis, there are two simple poles in the area bounded by
the contour, so by residue theorem:
\begin{align*}
B=&\frac{1}{2\pi
i}\int_{\gamma-\gamma_{\infty}}\frac{dE(q,\overline{q},y)}{\omega(q)}
\left(\zeta_n(y(q))B(y(\overline{q}),y_i)+\zeta_n(y(\overline{q}))B(y(q),y_i)\right)\\
=&\mathrm{Res}_{z(q)=t^{-1},s(t)^{-1}}
\frac{\frac{f+1}{2}(\frac{1}{z-z(q)}-\frac{1}{z-z(\overline{q})})}
{\log(1+\frac{z(q)}{f})-\log(1+\frac{z(\overline{q})}{f})}dydy_i\cdot
\frac{t(q)(t(q)-1)(ft(q)+1)}{(f+1)}\\
&\cdot\left(
F_n(t(q))\frac{dy(\overline{q})}{(y(\overline{q})-y_i)^2}+F_n(t(\overline{q}))
\frac{s'(t(q))dy(q)}{(y(q)-y_i)^2}\right)\\
=&\frac{-s'(t)dtdt_i}{2(\log(1+\frac{1}{ft})-\log(1+\frac{1}{fs(t)}))}
\cdot\frac{t(t-1)(ft+1)}{(f+1)}
\cdot\left(\frac{F_n(t)}{(s(t)-t_i)^2}+\frac{F_n(s(t))}{(t-t_i)^2}\right)\\
&+\frac{-dtdt_i}{2(\log(1+\frac{1}{ft})-\log(1+\frac{1}{fs(t)}))}
\cdot\frac{s(t)(s(t)-1)(fs(t)+1)}{(f+1)}
\cdot\left(\frac{F_n(s(t))}{(t-t_i)^2}+\frac{F_n(t)}{(s(t)-t_i)^2}\right)\\
=&\frac{-s'(t)dtdt_i}{\log(1+\frac{1}{ft})-\log(1+\frac{1}{fs(t)})}
\cdot\frac{t(t-1)(ft+1)}{(f+1)}
\cdot\left(\frac{F_n(t)}{(s(t)-t_i)^2}+\frac{F_n(s(t))}{(t-t_i)^2}\right)\\
=&\frac{-dtdt_i}{\log(1+\frac{1}{ft})-\log(1+\frac{1}{fs(t)})}
\left(\frac{\phi_{n+1}(t)s'(t)}{(s(t)-t_i)^2}+\frac{\phi_{n+1}(s(t))}{(t-t_i)^2}\right)\\
=&-\frac{\phi_{n+1}(t)B(s(t),t_i)+\phi_{n+1}(s(t))B(t,t_i)}
{\log(1+\frac{1}{ft})-\log(1+\frac{1}{fs(t)})}.
\end{align*}
After a more careful examination of the expression $B$, by expand it
as a formal Laurent series, it is not hard to see that in fact
$B\in\mathbb{C}(f)[t_i]((t))$. For a Laurent series
$$
f(t,t_i)=\sum_{n\geq
-N}^{\infty}a_n(t_i)t^n\in\mathbb{C}(f)[t_i]((t)),
$$
we denote by $f(t,t_i)_+$ the truncation $\sum_{n\geq
0}^{\infty}a_n(t_i)t^n\in\mathbb{C}(f)[t_i][[t]]$ of the regular
part of $f(t,t_i)$. By a similar argument as we have done for the
Type I residue, we have the estimate
$$
|\frac{1}{2\pi
i}\int_{\gamma}\frac{dE(q,\overline{q},y)}{w(q)}(\zeta_a(y(q))B(y(\overline{q}),y_i)
+\zeta_a(y(\overline{q}))B(y(q),y_i))|\leq
M\cdot\frac{dtdt_i}{|t^2|}\thicksim\frac{dtdt_i}{|t^2|},
$$
which has to be the principle part of the Laurent series, according
to Lemma \ref{Lemma2}. This finishes the proof the the following
\begin{proposition}
$$
P_n(t,t_i)=\left(\frac{\phi_{n+1}(t)B(s(t),t_i)+\phi_{n+1}(s(t))B(t,t_i)}
{\log(1+\frac{1}{ft})-\log(1+\frac{1}{fs(t)})}\right)_+
$$
\end{proposition}


\vskip 1cm

The calculation of this section proves the following theorem:
\begin{theorem} \label{Thm1}
The topological recursion in the BKMP conjecture \ref{Conj1} is
equivalent to
\begin{align*}
&(f(f+1))^{n-1}\sum_{b_1,\cdots,b_n}<\tau_{b_1}\cdots\tau_{b_n}\Gamma_g(f)>_g
\prod_{k=1}^{n}d\phi_{b_k}(t_k)\\
= &(f(f+1))^{n}\sum_{a_1,a_2,b_2,\cdots,b_n}
<\tau_{a_1}\tau_{a_2}\prod_{k=2}^{n}\Gamma_{g-1}(f)>_{g-1}
\prod_{k=2}^{n}d\phi_{b_k}(t_k)\cdot P_{a_1,a_2}(t_1)dt_1\\
&-(f(f+1))^{n-1}\sum_{g_1+g+2=g,I\coprod J=\{2,\cdots,n\}}^{stable}
\sum_{a_1,a_2,b_2,\cdots,b_n}<\tau_{a_1}\prod_{i\in
I}\tau_{b_i}\Gamma_{g_1}(f)>_{g_1}\\
&\cdot<\tau_{a_2}\prod_{j\in J}\tau_{b_j}\Gamma_{g_2}(f)>_{g_2}
\prod_{k=2}^{n}d\phi_{b_k}(t_k)\cdot P_{a_1,a_2}(t_1)dt_1\\
&-(f(f+1))^{n-2}\sum_{j=2}^n \sum_{b,b_i,i\in\{2,\cdots,n\}\setminus
\{j\}}<\tau_{b}\prod_{k=2,k\neq j}^n\tau_{b_k}>\prod_{k=2,k\neq
j}^{n}d\phi_{b_k}(t_k)\cdot P_b(t_1,t_j).
\end{align*}
\end{theorem}

\vskip 1.5cm

\section{Some Formal Analysis} \label{Sec5}

The critical point of the curve $x=y^f(1-f)$ is
$$
(x,y)=(\frac{f^f}{(f+1)^{f+1}},\frac{f}{f+1}). $$
Let
\begin{align*}
x=\frac{f^f}{(f+1)^{f+1}}e^{-\frac{f+1}{f}w} \hskip 1cm \text{and}
\hskip 1cm w=\frac{1}{2}v^2,
\end{align*}
then we have
$$
v^2=z^2\left(1+\sum_{k=1}^{+\infty}\frac{z^k}{k+2}(\frac{1-(\frac{-1}{f})^{k+1}}
{1-(\frac{-1}{f})})\right)\in\mathbb{Q}(f)[[z]].
$$

\vskip 0.3cm

Let $F(z)$ be the unique formal power series in $\mathbb{Q}(f)[[z]]$
such that $F(0)=1$ and
$$
F(z)^2=\left(1+\sum_{k=1}^{+\infty}\frac{z^k}{k+2}(\frac{1-(\frac{-1}{f})^{k+1}}
{1-(\frac{-1}{f})})\right).
$$
The coefficients of $F(z)F(z)=1+\sum_{k=1}^{\infty}a_kz^k$ are
determined recursively by the relations
\begin{align*}
&a_0=1\\
&\sum_{i=0}^ka_ia_{k-i}=\frac{1}{k+2}\left(\frac{1-(\frac{-1}{f})^{k+1}}{1+\frac{1}{f}}\right)
\hskip 1cm \text{for} \hskip 0.5cm k\geq 1.
\end{align*}

\vskip 0.3cm

It is obvious that the two possible solutions of $v$ are $\pm
zF(z)$. We take $v$ to be the solution $zF(z)$, which admit a formal
inverse function
$$
z=vG(v)=v+\sum_{k=1}^{\infty}g_k v^{k+1}\in v\mathbb{Q}(f)[[v]].
$$
The holomorphic map $s$ can be described as by sending
$t=\frac{1}{vG(v)}$ to $s(t)=\frac{1}{-vG(-v)}$, i.e, sending
$v\mapsto -v$. It is easy to see that
$$
t=\frac{1}{v[1+\sum_{k=1}^{\infty}g_kv^k]}=\frac{1}{v}[1+\sum_{k=1}^{\infty}h_kv^k]\in
v^{-1}\mathbb{Q}(f)[[v]].
$$

\vskip 0.5cm

\begin{lemma} \label{Lemma3}
Let
$$
\eta_{-1}(v)=-\frac{1}{2}(\log(1+\frac{1}{ft})-\log(1+\frac{1}{fs(t)}))
=\frac{1}{2}(\phi_{-1}(t)-\phi_{-1}(s(t))),
$$
and let
$\eta_{n+1}(v)=-\frac{f}{f+1}\frac{1}{v}\frac{d}{dv}\eta_n(v)$ for
all $n\geq -1$, then we have
$$
\eta_n(v)=\frac{1}{2}(\phi_n(t)-\phi_n(s(t))).
$$
Moreover, for each $n\geq -1$
$$
\eta_n(v)=\phi_n(t)+F_n(w)
$$
for some $F_n(w)\in\mathbb{Q}(f)[[w]]$.
\end{lemma}

\begin{proof}
By definition, we have
$$
-\frac{f}{f+1}\cdot\frac{1}{v}\frac{d}{dv}
=-\frac{f}{f+1}\cdot\frac{d}{dw}=x\frac{d}{dx}=\frac{t(t-1)(ft+1)}{(f+1)}\frac{d}{dt}
$$
is invariant under the action of $s:t\mapsto s(t)$. Applying this
operator $n$ times to the definition of $\eta_{-1}(v)$, we get the
first claim.

It is obvious that $\eta_n(v)$ is odd in the variable $v$, so
$\eta_{-1}(v)=\frac{-v}{f}(1+\cdots)\in v\mathbb{Q}(f)[[v]]$
contains only odd power of $v$. By applying the operator
$-\frac{f}{f+1}\frac{1}{v}\frac{d}{dv}$ to $\eta_{-1}(v)$, we get
$$
\eta_n(v)=\frac{f^n}{(f+1)^{n+1}}\cdot
\frac{(2n-1)!!}{v^{2n+1}}(1+\cdots)\in
\frac{1}{v^{2n+1}}\mathbb{Q}(f)[[w]].
$$
We have
\begin{align*}
\eta_{-1}(v)-\phi_{-1}(t)&=-\frac{1}{2}(\phi_{-1}(t)+\phi_{-1}(s(t)))\\
&=\frac{1}{2}\sum_{n=1}^{+\infty}\frac{(-1)^{n-1}}{nf^n}
\left(v^nG(v)^n+(-v)^nG(-v)^n\right),
\end{align*}
which is in the formal power series ring $\mathbb{Q}(f)[[v]]$.
Moreover, it is even in the variable $v$, thus in the ring
$\mathbb{Q}(f)[[w]]$. We denote it by $F_{-1}(w)$. Since there is no
constant term, $F_{-1}(w)\in \mathbb{Q}(f)[[w]]$. This proves the
second claim for $n=-1$.

For $n>-1$, apply the operator
$-\frac{f}{f+1}\frac{1}{v}\frac{d}{dv}$ to the equation
$$
\eta_{-1}(v)-\phi_{-1}(t)=F_{-1}(w)
$$
repeatedly $n$ times, we obtain
$$
\eta_n(v)=\phi_n(t)+F_n(w)
$$
for
$F_n(w)=(-\frac{f+1}{f}\frac{d}{dw})^nF_{-1}(w)\in\mathbb{Q}(f)[[w]]$
by induction.
\end{proof}

\vskip 0.5cm

\begin{corollary}
$$
P_{a,b}dt=-\frac{1}{2}\left(\frac{\eta_{a+1}(v)\eta_{b+1}(v)}{\eta_{-1}(v)}
(\frac{f+1}{f})vdv|_{v=v(t)}\right)_+
$$
\end{corollary}

\vskip 0.5cm

\begin{proof}
By Proposition \ref{Prop1}, we have
\begin{align*}
&\frac{\frac{1}{2}}{\log(1+\frac{1}{ft})-\log(1+\frac{1}{fs(t)})}
\left(\phi_{a+1}(t)\phi_{b+1}(s(t))+\phi_{a+1}(s(t))\phi_{b+1}(t)
\right)\cdot \frac{(f+1)dt}{t(t-1)(ft+1)}\\
&=\frac{1}{2\eta_{-1}(v)}\cdot \frac{(f+1)dt}{t(t-1)(ft+1)}
\cdot(\frac{\phi_{a+1}(t)-\phi_{a+1}(s(t))}{2}
\frac{\phi_{b+1}(t)-\phi_{b+1}(s(t))}{2}\\
&~~~-\frac{\phi_{a+1}(t)+\phi_{a+1}(s(t))}{2}
\frac{\phi_{b+1}(t)+\phi_{b+1}(s(t))}{2})\\
&=-\frac{(f+1)vdv}{2f\eta_{-1}(v)}
(\eta_{a+1}(v)\eta_{b+1}(v)-F_{a+1}(w)F_{b+1}(w)).
\end{align*}
To complete the proof, only need to notice that
$$
\left(\frac{F_{a+1}(w)F_{b+1}(w)vdv}{\eta_{-1}(v)}|_{v\mapsto
v(t)}\right)_+=0.
$$
\end{proof}

\vskip 0.5cm

\section{The left hand side} \label{Sec6}

After all the preparations, we are now ready to prove the BKMP
conjecture. We will show that the push forward of the Symmetrized
Cut-Join Equation via the first variable from the Framed Mirror
Curve to $\mathbb{C}$ gives the BKMP conjecture.

In this section, we will deal with the LHS of the Symmetrized
Cut-Join Equation. Unlike the Hurwitz number case studied in
\cite{EMS}, the LHS of our Symmetrized Cut-Join Equation involves
taking derivative with respect to the parameter $f$ that need more
special treatment.

The LHS of the symmetrized cut-join equation reads
\begin{align*}
&\left(\frac{\partial}{\partial
f}|_{t_i}+\sum_{l=1}^n\frac{t_l(t_l-1)}{f+1}\cdot\frac{\partial}{\partial
t_l}\right)H_n^g(t_1,\cdots,t_n,f)\\
&=\frac{\partial}{\partial f}|_{y_i}
H_n^g(t_1(y_1,f),\cdots,t_n(y_n,f),f)\\
&=\left(\frac{\partial}{\partial f}|_{x_i}+\sum_{l=1}^{n}\log
y_l\cdot x_l\frac{\partial}{\partial x_l} \right)H_n^g.
\end{align*}
If we change the variable
$$
x=\frac{f^f}{(f+1)^{f+1}}e^{-(\frac{f+1}{2f})v^2},
$$
then
$$
\frac{dx}{df}=x\left(\frac{f}{(f+1)^2}+\frac{v^2}{2f^2}+\log(\frac{f}{f+1})\right),
$$
and the differential operator above then becomes
\begin{align*}
&\frac{\partial}{\partial f}|_{x_i}+\sum_{l=1}^{m}\log y_l\cdot
x_l\frac{\partial}{\partial x_l}\\
&=\frac{\partial}{\partial f}|_{v_i}-\sum_{l=1}^m\left(
\log(1+\frac{1}{ft_l})-\frac{f}{(f+1)^2}-\frac{v_l^2}{2f^2}\right)
\frac{f}{v_l(f+1)}\frac{\partial}{\partial v_l}.
\end{align*}

\vskip 0.5cm

Taking the direct image of a function $F$ on the framed curve via
the map $\pi: C\rightarrow \mathbb{C}$ with respect to the first
variable, i.e, $\pi_*(F)=F(v_1)+F(-v_1)$. The term with
$\frac{\partial}{\partial f}$ becomes
\begin{align*}
&\frac{\partial}{\partial f}|_{v_i}\left(-(f(f+1))^{n-1}
\sum_{b_1,\cdots,b_n}<\tau_{b_1}\cdots\tau_{b_n}\Gamma_g(f)>_{g}
(\phi_{b_1}(t_1)+\phi_{b_1}(s(t_1)))\prod_{k=2}^n\phi_{b_k}(t_k)\right)\\
=&\frac{\partial}{\partial f}|_{v_i}\left(-2(f(f+1))^{n-1}
\sum_{b_1,\cdots,b_n}<\tau_{b_1}\cdots\tau_{b_n}\Gamma_g(f)>_{g}F_{b_1}(w_1)\prod_{k=2}^n\phi_{b_k}(t_k)
\right),
\end{align*}
which we denoted by $\widetilde{A}$ for short. In the summation, the
$l\neq 1$ term becomes a sum of
\begin{align*}
\left(\log(1+\frac{1}{ft_l})-\frac{f}{(f+1)^2}-\frac{v_l^2}{2f^2}\right)
\frac{f}{(f+1)v_l}\frac{\partial}{\partial
v_l}\left(F_{w_1}(t_1)\prod_{k=2}^n\phi_{b_k}(t_k)\right),
\end{align*}
which we denote by $A_l(b_1,\cdots,b_n)$. The term with $l=1$
contributes a sum of
\begin{align*}
&\left(\log(1+\frac{1}{ft_1})-\frac{f}{(f+1)^2}-\frac{v_1^2}{2f^2}\right)
\frac{f}{(f+1)v_1}\frac{\partial}{\partial
v_1}(\eta_{b_1}(v_1)-F_{b_1}(w_1))\prod_{k=2}^{n}\phi_{b_k}(t_k)\\
+&\left(\log(1+\frac{1}{fs(t_1)})-\frac{f}{(f+1)^2}-\frac{v_1^2}{2f^2}\right)
\frac{f}{(f+1)v_1}\frac{\partial}{\partial
v_1}(\eta_{b_1}(-v_1)-F_{b_1}(w_1))\prod_{k=2}^{n}\phi_{b_k}(t_k)\\
=&\left(\log(1+\frac{1}{ft_1})-\log(1+\frac{1}{fs(t_1)})\right)\frac{f}{(f+1)v_1}\frac{\partial}{\partial
v_1}\eta_{b_1}(v_1)\prod_{k=2}^{n}\phi_{b_k}(t_k)\\
+&\left(\log(1+\frac{1}{ft_1})+\log(1+\frac{1}{fs(t_1)})-\frac{2f}{(f+1)^2}-\frac{v_1^2}{f^2}\right)
F_{b_1+1}(w_1)\prod_{k=2}^{n}\phi_{b_k}(t_k)\\
=&2\eta_{-1}(v_1)\eta_{b_1+1}(v_1)\prod_{k=2}^{n}\phi_{b_k}(t_k)+A_1(b_1,\cdots,b_n),
\end{align*}
where $A_1(b_1,\cdots,b_n)$ denote the part involves
$F_{b_1+1}(w_1)$.

In summary, the push forward of the LHS of the Symmetrized Cut-Join
Equation is the following:
\begin{align*}
L:=&\widetilde{A}+(f(f+1))^{n-1}\sum_{b_1,\cdots,b_n}<\tau_{b_1}\cdots\tau_{b_n}\Gamma_g(f)>_g\\
&\cdot\left(2\eta_{-1}(v_1)\eta_{b_1+1}(v_1)\prod_{k=2}^{n}\phi_{b_k}(t_k)
+\sum_{l=1}^nA_l(b_1,\cdots,b_n)\right)\\
=&2(f(f+1))^{n-1}\sum_{b_1,\cdots,b_n}<\tau_{b_1}\cdots\tau_{b_n}\Gamma_g(f)>_g
\eta_{-1}(v_1)\eta_{b_1+1}(v_1)\prod_{k=2}^{n}\phi_{b_k}(t_k)+\widetilde{B},
\end{align*}
where we write
\begin{align*}
\widetilde{B}:=\widetilde{A}+(f(f+1))^{n-1}\sum_{b_1,\cdots,b_n}
<\tau_{b_1}\cdots\tau_{b_n}\Gamma_g(f)>_g\sum_{l=1}^nA_l(b_1,\cdots,b_n).
\end{align*}
for short. The crucial observation is that $\widetilde{A}$ and
$A_l(b_1,\cdots,b_n)$ for $l=1,\cdots,n$ are all regular in $w_1$.
This can be check from their definition. Thus their sum
$\widetilde{B}$ is also regular in $w_1$, and we have the important
estimate
$$
\left(\frac{\widetilde{B}}{2\eta_{-1}(v_1)}\cdot
\frac{(f+1)dt}{t(t-1)(ft+1)}\right)_+=0.
$$

\vskip 0.5cm

On the other hand, we have
\begin{align*}
W_g=&(-1)^{g+n}(f(f+1))^{n-1}\sum_{b_1,\cdots,b_n}
<\tau_{b_1}\cdots\tau_{b_n}\Gamma_g(f)>_g
\prod_{k=1}^{n}d\phi_{b_k}(t_k)\\
&=(-1)^{g+n}\frac{(f+1)dt_1}{2\eta_{-1}(v_1)t_1(t_1-1)(ft_1+1)}d_2\cdots
d_n \\
&\cdot\left(L-\widetilde{B}+(f(f+1))^{n-1}\sum_{b_1,\cdots,b_n}
<\tau_{b_1}\cdots\tau_{b_n}\Gamma_g(f)>_gF_{b_1+1}(w_1)\prod_{k=2}^{n}\phi_{b_k}(t_k)\right)\\
&=(-1)^{g+n}\left(\frac{(f+1)dt_1}{2\eta_{-1}(v_1)t_1(t_1-1)(ft_1+1)}d_2\cdots
d_n L\right)_+
\end{align*}
This finishes the computation of the LHS of the BKMP conjecture.

\vskip 2cm

\section{The Right Hand Side} \label{Sec7}


We now compute the push forward of the RHS of the Symmetrized
Cut-Join Equation. Recall that it is the sum of the following four
terms:
\begin{align*}
T_1=&-\frac{1}{2}\sum_{l=1}^nt_l(t_l-1)(\frac{t_lf+1}{f+1})\frac{\partial}{\partial
t_l}\cdot
t_{n+1}(t_{n+1}-1)(\frac{t_{n+1}f+1}{f+1})\frac{\partial}{\partial
t_{n+1}}H_{n+1}^{g-1}|_{t_{n+1}=t_l}\\
T_2=&-\frac{1}{2}\sum_{1\leqslant a\leqslant g-1}\sum_{1\leqslant
k\leqslant
n}\Theta_{k-1}(t_1(t_1-1)(\frac{t_1f+1}{f+1})\frac{\partial}{\partial
t_1}H_k^a(t_1,\cdots,t_k,f))\\
&\cdot(t_1(t_1-1)(\frac{t_1f+1}{f+1})\frac{\partial}{\partial
t_1}H_{n-k+1}^{g-a}(t_1,t_{k+1},\cdots,t_n,f))\\
T_3=&-\sum_{k=3}^n\Theta_{k-1}(t_1(t_1-1)(\frac{t_1f+1}{f+1})\frac{\partial}{\partial
t_1}H_k^0(t_1,\cdots,t_k,f))(t_1(t_1-1)
(\frac{t_1f+1}{f+1})\frac{\partial}{\partial
t_1}H_{n-k+1}^{g}(t_1,t_{k+1},\cdots,t_n,f))\\
T_4=&\Theta_1\frac{t_1^2(t_1-1)(t_2-1)}{t_1-t_2}(\frac{t_1f+1}{f+1})^2\frac{\partial}{\partial
t_1}H_{n-1}^{g}(t_1,t_3,\cdots,t_n,f).
\end{align*}

\vskip 0.4cm

The push forward of the $T_1$ part of the RHS of the Symmetrized
Cut-Join Equation is the sum of
\begin{align*}
\frac{1}{2}(f(f+1))^{n}&\sum_{b_1,\cdots,b_{n+1}}
<\tau_{b_1}\cdots\tau_{b_{n+1}}\Gamma_{g-1}(f)>_{g-1}\\
&\cdot(\phi_{b_1+1}(t_1)\phi_{b_{n+1}+1}(t_1)+\phi_{b_1+1}(s(t_1))\phi_{b_{n+1}+1}(s(t_1))
)\cdot \prod_{k=2}^n\phi_{b_k}(t_k)\\
=(f(f+1))^{n}&\sum_{b_1,\cdots,b_{n+1}}
<\tau_{b_1}\cdots\tau_{b_{n+1}}\Gamma_{g-1}(f)>_{g-1}\\
&\cdot(\eta_{b_1+1}(v_1)\eta_{b_{n+1}+1}(v_1)+F_{b_1+1}(w_1)F_{b_{n+1}+1}(w_1))
\cdot\prod_{k=2}^n\phi_{b_k}(t_k)
\end{align*}
and
\begin{align*}
-\sum_{l=2}^m(f(f+1))^{n}\sum_{b_1,\cdots,b_{n+1}}
<\tau_{b_1}\cdots\tau_{b_{n+1}}\Gamma_{g-1}(f)>_{g-1}
F_{b_1}(w_1)(\phi_{b_l+1}(t_l)\phi_{b_{n+1}+1}(t_l))
\prod_{i=2,i\neq l}^{n}\phi_{b_i}(t_i).
\end{align*}
Thus we have
\begin{align*}
&\left(\frac{\pi_* T_1\cdot
(f+1)dt_1}{2\eta_{-1}(v_1)t_1(t_1-1)(ft_1+1)}\right)_+\\
=&(f(f+1))^{n}\sum_{b_1,\cdots,b_{n+1}}
<\tau_{b_1}\cdots\tau_{b_{n+1}}\Gamma_{g-1}(f)>_{g-1}\prod_{k=2}^n\phi_{b_k}(t_k)\\
&\cdot\left(-\frac{(f+1)\eta_{b_1+1}(v_1)\eta_{b_{n+1}+1}(v_1)v_1dv_1}{2f\eta_{-1}(v_1)}
\right)_+
\\
=&(f(f+1))^{n}\sum_{b_1,\cdots,b_{n+1}}
<\tau_{b_1}\cdots\tau_{b_{n+1}}\Gamma_{g-1}(f)>_{g-1}
\prod_{k=2}^n\phi_{b_k}(t_k)\cdot P_{b_1,b_{n+1}}(t_1)dt_1.
\end{align*}

\vskip 0.5cm

By completely the same computation, applied to $T_2+T_3$, we have
\begin{align*}
&\left(\frac{\pi_* (T_2+T_3)\cdot
(f+1)dt_1}{2\eta_{-1}(v_1)t_1(t_1-1)(ft_1+1)}\right)_+\\
=&-(f(f+1))^{n-1}\sum_{g_1+g_2=g} \sum_{I\coprod
J=\{2,\cdots,n\}}^{stable}\sum_{a_1,a_2,b_2,\cdots,b_n}
<\tau_{a_1}\prod_{i\in
I}\tau_{b_i}\Gamma_{g_1}(f)>_{g_1}\\
&\cdot <\tau_{a_2}\prod_{j\in
J}\tau_{b_j}\Gamma_{g_2}(f)>_{g_2}P_{a_1,a_2}(t_1)dt_1
\prod_{k=2}^{n}\phi_{b_k}(t_k).
\end{align*}

\vskip 0.8cm

The $T_4$ part is much more complicated, and need additional care.
First of all, it is equal to
\begin{align*}
&-\sum_{1\leq i,j\leq n}\frac{t_i(ft_i+1)(t_j-1)}{(f+1)(t_i-t_j)}
(f(f+1))^{n-2}\sum_{b,b_i,i\in\{1,\cdots,n\}\setminus
\{i,j\}}<\tau_b\prod_{k=1,k\neq i,j}^{n}\tau_{b_k}\Gamma_{g}(f)>_{g}
\phi_{b+1}(t_i)\prod_{k=1,k\neq i,j}^{n}\phi_{b_k}(t_k).\\
\end{align*}
The push forward $\pi_* T_4=S_1+S_2+S_3$, for
\begin{align*}
S_1=&\sum_{2\leq i\neq j\leq
n}2(f(f+1))^{n-2}\frac{t_i(ft_i+1)(t_j-1)}{(f+1)(t_i-t_j)}
\sum_{b,b_i,i\in\{1,\cdots,n\}\setminus
\{i,j\}}<\tau_b\prod_{k=1,k\neq
i,j}^{n}\tau_{b_k}\Gamma_{g}(f)>_{g}\\
&\cdot F_{b_1}(w_1)\phi_{b+1}(t_i)\prod_{k=2,k\neq i,j}^{n}\phi_{b_k}(t_k)\\
S_2=&-\sum_{j=2}^n(f(f+1))^{n-2}\sum_{b,b_i,i\in\{2,\cdots,n\}\setminus
\{j\}}<\tau_b\prod_{k=2,k\neq j}^{n}\tau_{b_k}\Gamma_{g}(f)>_{g}\\
&\cdot\left(\frac{t_1(ft_1+1)(t_j-1)}{(f+1)(t_1-t_j)}\phi_{b+1}(t_1)
+\frac{s(t_1)(fs(t_1)+1)(t_j-1)}{(f+1)(s(t_1)-t_j)}\phi_{b+1}(s(t_1))\right)
\prod_{k=2,k\neq j}^{n}\phi_{b_k}(t_k)\\
S_3=&-\sum_{i=2}^n(f(f+1))^{n-2}\sum_{b,b_i,i\in\{2,\cdots,n\}\setminus
\{i\}}<\tau_b\prod_{k=2,k\neq i}^{n}\tau_{b_k}\Gamma_{g}(f)>_{g}\\
&\cdot\left(\frac{t_i(ft_i+1)(t_1-1)}{(f+1)(t_i-t_1)}
+\frac{t_i(ft_i+1)(s(t_1)-1)}{(f+1)(t_i-s(t_1))}\right)
\phi_{b+1}(t_i)\prod_{k=2,k\neq i}^{n}\phi_{b_k}(t_k).
\end{align*}

It is easy to see that
$$
\left(\frac{S_1\cdot
(f+1)dt_1}{2\eta_{-1}(v_1)t_1(t_1-1)(ft_1+1)}\right)_+=0,
$$
since $S_1$ is holomorphic in the variable $w_1$. Because of
\begin{align*}
&\frac{t_i(ft_i+1)(t_1-1)}{(f+1)(t_i-t_1)}
+\frac{t_i(ft_i+1)(s(t_1)-1)}{(f+1)(t_i-s(t_1))}\\
=&-2\frac{t_i(ft_i+1)}{f+1}-\frac{t_i(t_i-1)(ft_i)}{f+1}
\cdot\sum_{k=0}^{\infty}\left(\frac{t_i^k}{t_1^{k+1}}+\frac{t_i^k}{s(t_1)^{k+1}}\right),
\end{align*}
we also have
$$
\left(\frac{S_3\cdot
(f+1)dt_1}{2\eta_{-1}(v_1)t_1(t_1-1)(ft_1+1)}\right)_+=0.
$$

\vskip 0.5cm

Finally, let us compute the contribution of $S_2$. Consider the
differentials:
\begin{align*}
&d_j\left(\frac{(f+1)dt_1}{2\eta_{-1}(v_1)t_1(t_1-1)(ft_1+1)}\cdot
\frac{t_1(ft_1+1)(t_j-1)}{(f+1)(t_1-t_j)}\phi_{b+1}(t_1)\right)\\
=&d_j\left(\frac{(t_j-1)dt_1}{2\eta_{-1}(v_1)(t_1-1)(t_1-t_j)}\phi_{b+1}(t_1)\right)\\
=&d_j\left(-\frac{\phi_{b+1}(t_1)dt_1}{2\eta_{-1}(v_1)(t_1-1)}
+\frac{\phi_{b+1}(t_1)dt_1}{2\eta_{-1}(v_1)(t_1-t_j)}\right)\\
=&\frac{\phi_{b+1}(t_1)B(t_1,t_j)}{2\eta_{-1}(v_1)},\\
&d_j\left(\frac{(f+1)dt_1}{2\eta_{-1}(v_1)t_1(t_1-1)(ft_1+1)}\cdot
\frac{s(t_1)(fs(t_1)+1)(t_j-1)}{(f+1)(s(t_1)-t_j)}\phi_{b+1}(s(t_1))\right)\\
=&d_j\left(\frac{(f+1)s'(t_1)dt_1}{2\eta_{-1}(v_1)s(t_1)(s(t_1)-1)(fs(t_1)+1)}\cdot
\frac{s(t_1)(fs(t_1)+1)(t_j-1)}{(f+1)(s(t_1)-t_j)}\phi_{b+1}(s(t_1))\right)\\
=&d_j\left(\frac{(t_j-1)s'(t_1)dt_1}{2\eta_{-1}(v_1)(s(t_1)-1)(s(t_1)-t_j)}\phi_{b+1}(s(t_1))\right)\\
=&d_j\left(-\frac{\phi_{b+1}(s(t_1))s'(t_1)dt_1}{2\eta_{-1}(v_1)(s(t_1)-1)}
+\frac{\phi_{b+1}(s(t_1))s'(t_1)dt_1}{2\eta_{-1}(v_1)(s(t_1)-t_j)}\right)\\
=&\frac{\phi_{b+1}(s(t_1))B(s(t_1),t_j)}{2\eta_{-1}(v_1)}.
\end{align*}
We have the following simplified expression
\begin{align*}
&d_2\cdots d_n\left(\frac{S_2\cdot
(f+1)dt_1}{2\eta_{-1}(v_1)t_1(t_1-1)(ft_1+1)}\right)\\
=&-\sum_{j=2}^n(f(f+1))^{n-2}\sum_{b,b_i,i\in\{2,\cdots,n\}\setminus
\{j\}}<\tau_b\prod_{k=2,k\neq
j}^{n}\tau_{b_k}\Gamma_{g}(f)>_{g}\prod_{k=2,k\neq j}^{n}d\phi_{b_k}(t_k)\\
&\cdot\left(\frac{\phi_{b+1}(t_1)B(t_1,t_j)+\phi_{b+1}(s(t_1))B(s(t_1),t_j)}{2\eta_{-1}(v_1)}\right).\\
\end{align*}
To obtain the $P_n(t_1,t_j)$ term in the BKMP conjecture, we need to
switch one of the $t_1$ with $s(t_1)$ in the above expression. This
can be done by the following estimate:
\begin{align*}
\left(\frac{B(t_1,t_j)}{2\eta_{-1}(v_1)}(\phi_{b+1}(t_1)+\phi_{b+1}(s(t_1)))\right)_+
=\left(-\frac{B(t_1,t_j)F_{b+1}(w_1)}{\eta_{-1}(v_1)}\right)_+=0,
\end{align*}
and similarly
\begin{align*}
\left(\frac{B(s(t_1),t_j)}{2\eta_{-1}(v_1)}(\phi_{b+1}(t_1)+\phi_{b+1}(s(t_1)))\right)_+
=\left(-\frac{B(s(t_1),t_j)F_{b+1}(w_1)}{\eta_{-1}(v_1)}\right)_+=0.
\end{align*}
Thus we have the relation
\begin{align*}
P_{b}(t_1,t_j)dt_1dt_j
=&\left(\frac{\phi_{b+1}(s(t_1))B(t_1,t_j)+\phi_{b+1}(t_1)B(s(t_1),t_j)}
{\log(1+\frac{1}{ft})-\log(1+\frac{1}{fs(t)})} \right)_+\\
=&\left(\frac{\phi_{b+1}(t_1)B(t_1,t_j)+\phi_{b+1}(s(t_1))B(s(t_1),t_j)}{2\eta_{-1}(v_1)}
\right)_+\\
\end{align*}
and the estimate of $\pi_*T_4$ and $S_2$:
\begin{align*}
&d_2\cdots d_n\left(\frac{\pi_*T_4\cdot
(f+1)dt_1}{2\eta_{-1}(v_1)t_1(t_1)(ft_1+1)}\right)_+=d_2\cdots
d_n\left(\frac{S_2\cdot
(f+1)dt_1}{2\eta_{-1}(v_1)t_1(t_1)(ft_1+1)}\right)_+\\
=&-(f(f+1))^{n-2}\sum_{j=2}^n\sum_{b,b_i,i\in\{2,\cdots,n\}\setminus
\{j\}}<\tau_b\prod_{k=1,k\neq
i,j}^{n}\tau_{b_k}\Gamma_{g}(f)>_{g}\cdot
P_{b}(t_1,t_j)\prod_{k=2,k\neq j}^nd\phi_{b_k}(t_k).
\end{align*}

\vskip 0.5cm

Collect all the pieces, we have finish the computation of the push
forward of the RHS of the Symmetrized Cut-Join Equation:
\begin{align*}
&d_2\cdots d_n\left(\frac{\pi_*(T_1+T_2+T_3+T_4)\cdot
(f+1)dt_1}{2\eta_{-1}(v_1)t_1(t_1)(ft_1+1)}\right)_+\\
= &(f(f+1))^{n}\sum_{a_1,a_2,b_2,\cdots,b_n}
<\tau_{a_1}\tau_{a_2}\prod_{k=2}^{n}\Gamma_{g-1}(f)>_{g-1}
\prod_{k=2}^{n}d\phi_{b_k}(t_k)\cdot P_{a_1,a_2}(t_1)dt_1\\
&-(f(f+1))^{n-1}\sum_{g_1+g+2=g,I\coprod J=\{2,\cdots,n\}}^{stable}
\sum_{a_1,a_2,b_2,\cdots,b_n}<\tau_{a_1}\prod_{i\in
I}\tau_{b_i}\Gamma_{g_1}(f)>_{g_1}\\
&\cdot<\tau_{a_2}\prod_{j\in J}\tau_{b_j}\Gamma_{g_2}(f)>_{g_2}
\prod_{k=2}^{n}d\phi_{b_k}(t_k)\cdot P_{a_1,a_2}(t_1)dt_1\\
&-(f(f+1))^{n-2}\sum_{j=2}^n \sum_{b,b_i,i\in\{2,\cdots,n\}\setminus
\{j\}}<\tau_{b}\prod_{k=2,k\neq j}^n\tau_{b_k}>\prod_{k=2,k\neq
j}^{n}d\phi_{b_k}(t_k)\cdot P_b(t_1,t_j).
\end{align*}
The recursion in Theorem \ref{Thm1} then follows from plugging in
the results of this and the previous section into the equality $
L=\pi_*(T_1+T_2+T_3+T_4)$.

\end{document}